\documentclass[9pt]{article}
\usepackage{amsmath}
\usepackage{latexsym}
\usepackage{amsfonts}
\usepackage{amssymb}

\def\la{\langle}
\def\ra{\rangle}

\renewcommand{\phi}{\varphi}

\newcommand{\nf}{\hspace*{-5pt}}

\def \beq { \begin{equation} }
\def \eeq { \end{equation} }

\def\overset#1\to#2{\mathrel{\mathop{#2}\limits^{#1}}}
\def\underset#1\to#2{\mathrel{\mathop{#2}\limits_{#1}}}

\def \varinjlim {\underset{\longrightarrow}\to{lim}}
\def \varprojlim {\underset{\longleftarrow}\to{lim}}

\def \li {\varinjlim}
\def \lp {\varprojlim}

\def \lli#1 {\mathrel{\mathop{\li}\limits_{#1}}}
\def \llp#1 {\mathrel{\mathop{\lp}\limits_{#1}}}

\def \empty {\emptyset}
\def \leqs {\leqslant}
\def \rest {\restriction}
\def \la {\leftarrow}

\def \rsa {\rightsquigarrow }
\def \lrsa {\leftrightsquigarrow }
\def \rat {\rightarrowtail }
\def \thra {\twoheadrightarrow }

\def \hookr {\hookrightarrow }

\newcommand{\hem}{\hspace*{1em}}

\newcommand{\hfl}{\hspace*{\fill}}

\newlength{\dede}     

\newcommand{\hp}{\hspace*{\parindent}} 
\newcommand{\nhp}{\hspace*{-\parindent}}

\newcommand{\N}{\mbox{$ \mathbb N$}}   

\newcommand{\sub}{\mbox{$\subseteq$}}

\def \empty {\emptyset}

\newcommand{\cK}{\mbox{$\cal K$}}

\newcommand{\cM}{\mbox{$\cal M$}}

\newcommand{\all}{\mbox{$\forall$}}

\newcommand{\Ra}{\mbox{$\Rightarrow$}}
\newcommand{\Lra}{\mbox{$\Leftrightarrow$}}

\def \rest {\restriction}
\def \ra {\rightarrow}
\def \la {\leftarrow}

\def \rsa {\rightsquigarrow }
\def \lrsa {\leftrightsquigarrow }
\def \rat {\rightarrowtail }
\def \thra {\twoheadrightarrow }

\def \hookr {\hookrightarrow }

\def\overset#1\to#2{\mathrel{\mathop{#2}\limits^{#1}}}
\def\underset#1\to#2{\mathrel{\mathop{#2}\limits_{#1}}}

\newcommand{\und}[1]{\raisebox{-.2ex}{\underline{\raisebox{.2ex}{#1}}}}  

\def \leqs {\leqslant}
\def \geqs {\geqslant}

\def \varinjlim {\underset{\longrightarrow}\to{lim}}
\def \varprojlim {\underset{\longleftarrow}\to{lim}}

\def \li {\varinjlim}
\def \lp {\varprojlim}

\def \lli#1 {\mathrel{\mathop{\li}\limits_{#1}}}
\def \llp#1 {\mathrel{\mathop{\lp}\limits_{#1}}}

\newcommand{\td}[1]{\mbox{$\widetilde{#1}$}}

\newtheorem{Th}{Theorem}[section] 
\newtheorem{Co}[Th]{Corollary}
\newtheorem{Df}[Th]{Definition}
\newtheorem{Pro}[Th]{Proposition}
\newtheorem{Le}[Th]{Lemma}
\newtheorem{Exa}[Th]{Example}
\newtheorem{Rem}[Th]{Remark}
\newtheorem{Fa}[Th]{Fact}
\newtheorem{Que}[Th]{Question}
\newtheorem{Ct}[Th]{}
\newtheorem{Afi}[Th]{Claim}

\newcommand{\baf}{\begin{Afi}\nf{\sl . }}
\newcommand{\eaf}{\end{Afi}}

\newcommand{\bdf}{\begin{Df}\nf{\bf .}}
\newcommand{\edf}{\end{Df}}
\newcommand{\bte}{\begin{Th}\nf{\bf .}}
\newcommand{\ete}{\end{Th}}  
\newcommand{\bco}{\begin{Co}\nf{\bf .}}
\newcommand{\eco}{\end{Co}}
\newcommand{\ble}{\begin{Le}\nf{\bf .}}
\newcommand{\ele}{\end{Le}}
\newcommand{\bpr}{\begin{Pro}\nf{\bf .}}
\newcommand{\epr}{\end{Pro}}
\newcommand{\bex}{\begin{Exa}\nf{\bf .} \rm}
\newcommand{\eex}{\end{Exa}}
\newcommand{\bre}{\begin{Rem}\nf{\bf .} \rm}
\newcommand{\ere}{\end{Rem}}
\newcommand{\bfa}{\begin{Fa}\nf{\bf .} \sl}
\newcommand{\efa}{\end{Fa}}
\newcommand{\bqt}{\begin{Que}\nf{\bf .}}
\newcommand{\eqt}{\end{Que}}
\newcommand{\bdm}{\nhp{\bf Proof.  }}

\newcommand{\qdr}{\hfl $\square$ }
\newcommand{\bct}{\begin{Ct}\nf{\bf .} \rm}
\newcommand{\ect}{\end{Ct}}
\newcommand{\bxa}{\begin{Exa}\nf{\bf .} \rm}
\newcommand{\exa}{\end{Exa}}

\newcommand{\bu}{$\bullet$\ }

\begin{document}

\title{A global approach to AECs} 

\author{Hugo L. Mariano \thanks{Research supported by FAPESP, under the Thematic Project LOGCONS:
Logical consequence, reasoning and computation (number 2010/51038-0).} \thanks{Instituto de Matem\'atica e Estat\'istica, University of S\~ao Paulo, Brazil, email: hugomar@ime.usp.br }
  \hem
Pedro H. Zambrano  \hem \thanks{Departamento de Matemáticas, Faculdad de Ciencias, National University of Colombia, Colombia, emails: phzambranor@unal.edu.co, avillavecesn@unal.edu.co} \hem Andr\'es Villaveces}

\date{December 2013}

\maketitle

\begin{abstract}

In this work we present some general categorial ideas on Abstract Elementary Classes (AECs)
, inspired by the totality of AECs of the form $(Mod(T), \preceq)$ , for a first-order theory T: (i) we define a natural notion of (funtorial) morphism between AECs; (ii) explore the following constructions of AECs: "generalized" theories, pullbacks of AECs, (Galois) types as AECs; (iii) apply categorial and topological ideas to encode model-theoretic notions on spaces of types 
; (iv) present the "local" axiom for AECs here called  "local Robinson's property" and an application (Robinson's diagram method); (v) introduce the category $AEC$ of Grothendieck's gluings of all AECs (with change of basis); (vi) introduce the "global" axioms of "tranversal Robinson's property" (TRP) and "global Robinson's property" (GRP) and prove that TRP is equivalent to GRP and GRP entails a natural version of Craig interpolation property. 

\end{abstract}

\section{Introduction}








\hp The notion of Abstract Elementary Class (AEC) was introduced by Shelah in 1987 (\cite{She}) and is a an abstract model-theoretic setting (i.e., there is no syntax and proof theory associated to it)   that generalizes the natural occurrences:\\
 (i) $T \sub L_{\omega,\omega}-sent$\ $\rsa$\ $(Mod(T), \preceq)$;\hfl
 (ii) $\theta \in L_{\omega_1,\omega}-sent$\ $\rsa$\ $(Mod(\theta), \preceq)$\hfl

\bdf An AEC is a pair $(\cK, \preceq_K)$ where $\cK \sub L-str$, for some (finitary) language $L$ and $\preceq_K$ is a binary relation on $\cK$ that satisfies:

(i) \ {\em order}: \\
$\bullet$ $\preceq_K$ is a partial order relation on $\cK$;\hfl \hfl
$\bullet$ $M \preceq_K N$ $\Ra$ $M \sub_{L-str} N$.  

(ii) \ {\em isomorphism}: \\
$\bullet$ \ $M \in \cK$ and $M' \cong_{L-str} M$, then $M' \in \cK$; \\
$\bullet$ \ $M \preceq_K N$, $f\in L-iso(N,N')$ and $f\rest \in L-iso(M,M')$  
$\Ra$ $M' \preceq_K N'$ 

(iii) \ {\em coherence}:\\
$M, N, P \in \cK$, $M \sub_{L-str} N \sub_{L-str} P$ and $M \preceq_K P$, $N \preceq_K P$ $\Ra$ $M \preceq_K N$

(iv) \ {\em reunion}: Let  $\{ M_\alpha : \alpha \in \theta\} \sub \cK$, $M_\alpha \preceq_K M_\beta$ if $\alpha <  \beta < \theta$, then:\\
$\bullet$ \ $M := \bigcup_{\alpha \in \theta} M_\alpha \in \cK$ and $M_\alpha \preceq_K M$, $\all \alpha \in \theta$; \hfl  \hfl
$\bullet$ \ if $N \in \cK$ and $M_\alpha \preceq_K N$, $\all \alpha \in \theta$, then $M \preceq_K N$.  

(v) \ {\em L\"owenhiem-Skolem}:\
There is an infinite cardinal $\kappa$ such that for each $N \in \cK$ and $X \sub |N|$, there is $M \in \cK$, $M \preceq_K N$, $X \sub |M|$ and $\sharp |M| \leqs \sharp X + \kappa$. \
$LS(\cK) := min\{\kappa$: $\kappa$ satisfies the property$\}$
\qdr\edf

\bct Let  $(\cK, \preceq_K)$ be an AEC with $\cK \sub L-str$, then there are two natural notions:\\
\bu $\cK$-substructure: $C \preceq_K D$ $\Ra$ $C \overset{\preceq_K}\to\hookr D$ \hfl \hfl
\bu $\cK$-embedding: $(A \overset{f}\to\rat B)$ = $(A \overset{\cong_L}\to\ra A' \overset{\preceq_K}\to\hookr B)$ 
\qdr\ect


\bct {\em An AEC as a category:} Let $cat(\cK)$ be the category with $\cK$ as the class of objects and $\cK$-embeddings as arrows. 

\bu This is  well-defined:\\
\und{identities}: because $\preceq_K$ is reflexive\\
\und{composition}: by the closure under isos\\
\hfl $A \overset{\cong}\to\ra A' \overset{\preceq}\to\hookr B \overset{\cong}\to\ra B' \overset{\preceq}\to\hookr C$ \hfl\\
\hfl $A' \overset{\cong}\to\ra A'' \overset{\preceq}\to\hookr B'$ \hfl

\bu $cat(\cK)$ has all (upward) directed colimits.

\bu $cat(\cK)$ is an accessible category: \cite{AR}, \cite{Lie1} \cite{Lie2}, \cite{BR}.

\bu We will denote $cat(\cK)_{strict}$ the category with $\cK$ as the class of objects and $\preceq_K$-inclusions as arrows.
\qdr\ect





\bct {\em Shelah's Presentation Theorem} (\cite{She},  \cite{Bal2}):   Is the first occurence of "change of languages" in AECs-theory:\\
$\cK \sub L-str$ \hem $\Ra$ \hem
$\exists L' \supseteq L$  $\exists T' \sub L'-sent$ $\exists$ a set $\Gamma'$ of (finitary) $L'$-types such that,
if $EC(T', \Gamma') = \{ M' \in L'-str: M' \vDash T'$ and $ M'$ omits $\Gamma'\}$, then:\\
-- \ $\cK = \{M'\!\rest_L \in L-str: M' \in EC(T', \Gamma')\}$ \ (i.e. $\cK = PC_L(T',\Gamma')$);\\
-- \ $M' \preceq_{L'_{\omega\omega}} N' \in EC(T', \Gamma')$ \ $\Ra$\ $M'\!\rest_L \ \preceq_K N'\!\rest_L \in \cK$ 
\qdr\ect

\bct {\em Main Question: Shelah's categoricity conjecture} (\cite{Bal1}): \
 There is a cardinal $\mu(\kappa)$
such that for all AEC $\cK$ with  $LS(\cK) \leqs \kappa$, if $\cK$ is categorical in
some cardinal greater than $\mu(\kappa)$ then K is categorical in all $\lambda \geqs \mu(\kappa)$.
\qdr\ect

\bct {\em AECs and infinitary languages:} $LS(\cK) \leqs \kappa$ \hem $\Ra$ \hem $\cK$ is closed under $\equiv$ in $L_{\infty, \kappa^+}$ (\cite{Kue})
\qdr\ect


\bct {\em (Generalized) types in AECs:} The notion of type in AECs-theory (\cite{Lie1}, \cite{Lie2}) is given by another application  of "change of languages". We consider here an expansion of this notion, called "generalized type" or simply g-type (see, for instance, section I.4 in \cite{Van} for motivations).

{\em g-types}: $M \in \cK$, $I$ set, a (generalized) $I$-type over $M$:\\
is an equivalence class $(N,f,\bar{a})$ where $N \in \cK$, $f : M \rat N \in cat(\cK)$,  $\bar{a}: I \ra N$;\\
$(N,f,\bar{a}) \sim (N',f',\bar{a}')$ iff $\exists n \in \N$, $P_0, \ldots, P_n \in \cK$ and $\exists$ zigzag\\
 $(N,f,\bar{a}) = (P_0,f_0,\bar{a}_0) \overset{h_1}\to\ra (P_1,f_1,\bar{a}_1) \overset{h_2}\to\la (P_2,f_2,\bar{a}_2) \overset{h_3}\to\ra  \ldots \overset{h_{n-1}}\to\la (P_{n-1},f_{n-1},\bar{a}_{n-1}) \overset{h_{n}}\to\ra (P_n,f_n,\bar{a}_n) = (N',f',\bar{a}')$  \\  
$\bullet$ \ LS-axiom  $\Ra$ $\exists$ a set of representatives for $I$-types over $M$;\\ 
$\bullet$ \ if $cat(\cK)$ enjoys (AP), then $(N,f,\bar{a}) \sim (N',f',\bar{a}')$ iff $\exists P \in \cK$ and morfisms $h_1, h_2$ such that:\\
 $(N,f,\bar{a}) \overset{h_1}\to\ra (P,g,\bar{b}) \overset{h_2}\to\la (N',f',\bar{a}')$\\
$\bullet$ \ if $cat(\cK)$ enjoys (AP), (JEP), has arbitrary large models, then:\\
-- $\exists$ monster models;\\
-- $I$-types over $M$ $\lrsa$ orbits of $I$-sequences over a $\cK$-monster model over $M\cup I$ \ (Galois g-types).\\
notation: $g-T(M,I)$\\
\bu {\em topology on g-types}: this is analogous to the topology on usual types considered in \cite{Lie1} \\
$\lambda \geqs LS(\cK)$, $M \in \cK$ :\\
 $\cK^\lambda(M) := \{ N \in \cK: N \preceq_K M, \sharp N \leqs \lambda\}$,\\
basis of $X^\lambda(M,I)$, the $\lambda$-topology of  $g-T(M,I)$:\hem $\{ U_{p,N}$ :\ $p \in g-T(N,I)$, $N \in \cK^\lambda(M)\}$, where:\\
$U_{p,N} := \{ q \in g-T(M,I) : q\!\!\rest =p\}$\\
(remark: if $\sharp M \leqs \lambda$, then $X^\lambda(M,I)$ is discrete)
\qdr\ect

\section{The category of all AECs}

\bct {\em Motivation}: change of languages in model theory of FOL ({\bf Institution theory})\\
$\alpha: L \ra L'$ (morphism of languages) \hem $\rsa$ \hem $(\hat{\alpha}, \alpha^\star)$, where: \\
\bu $\hat{\alpha} : Form(L) \ra Form(L')$, a function defined by recursion on complexity\\
\bu  $\alpha^\star : L'-str \ra L-str$, a functor \\
 $|M'^\alpha| = |M'|, \ R^{M'^\alpha} = \alpha(R)^{M'}, \ f^{M'^\alpha} = \alpha(f)^{M'}, \  c^{M'^\alpha} = \alpha(c)^{M'}$\\
\bu   $M'^\alpha \vDash_L \theta [\bar{a}]$ $\Lra$ $M' \vDash_{L'} \hat{\alpha}(\theta)[\bar{a}]$\\
\bu  $f' : M' \ra N'$ is $L'$-homomorfism (respect. $L'$-embbeding, $L'$-elementary embedding) \hem $\Ra$ \\ $f': M'^\alpha \ra N'^\alpha$ is $L'$-homomorfism (respect. $L'$-embbeding, $L'$-elementary embedding)\\
\bu $T \sub L-sent$, $T' \sub L'-sent$ are closed theories such that $\hat{\alpha}^{-1}[T'] \supseteq T$, then $\alpha^\star\!\!\rest : (Mod(T'), \preceq_{L'}) \ra (Mod(T), \preceq_L)$ 
\qdr\ect

\bct {\em The category of all AECs} (For foundational issues $\rsa$ Grothendieck's universes)\\
\und{objects}: the "class of all" AECs\\
\und{arrows}: $\cK \sub L-str$, $\cK' \sub L'-str$ \ (remark that the language is determined from the (non-empty) class, set $\empty \sub \empty-str$)\\
$\alpha : L \ra L'$ such that the restriction $\alpha^\star\!\!\rest : (\cK',\preceq_{K'}) \ra (\cK, \preceq_{K})$ is well defined:

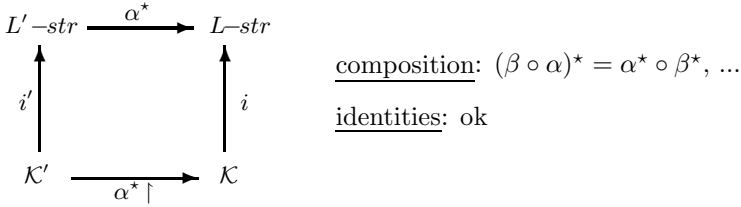
\begin{picture}(100,80)
\setlength{\unitlength}{.4\unitlength} \thicklines 
\put(-13,162){\small $L'\!-\!\!str$}
\put(5,20){\small $\cK'$} 
\put(180,162){\small $L\!\!-\!\!str$} 
\put(0,90){\small $i'$}
\put(100,175){\small $\alpha^\star$} 
\put(190,20){\small  $\cK$} 
\put(210,90){\small $i$} 
\put(90,05){\small $\alpha^\star\!\rest$}
\put(20,50){\vector(0,1){100}}
\put(195,50){\vector(0,1){100}} 
\put(65,167){\vector(1,0){100}} 
\put(50,25){\vector(1,0){120}}

\put(300,125){\und{composition}: $(\beta\circ\alpha)^\star = \alpha^\star \circ \beta^\star$, ...}

\put(300,75){\und{identities}: ok}
\end{picture}

{\em Example:} $inclusion : (K, \preceq_K) \hookr (L-str, \sub_L)$, $inclusion = (id_L)^\star\!\rest$.\\
\hp {\em Remark:} $\alpha^\star$ is always a faitful functor. If $\alpha$ is a "surjective" morphism of languages, then $\alpha^\star$ is a full functor.
\qdr\ect

\subsection{Constructions of AECs} 

\hp In this subsection we provide new examples of AECs, obtained from construction, for instance: elementary class (Mod(T)) in AECs; pullbacks of AECs; g-types in AECs, etc.

{\bf Theories Mod(T) in AEC}

\bct\ ({\bf FOL}) The $L$-satisfaction relation $\vDash_L$:\\
\bu establish a Galois connection $(Mod,Th)$ between sets of (first-order) $L$-sentences and classes of $L$-structures;\\
\bu $(Mod\!\rest,Th\!\rest)$ establish (order reversing) bijection between closed $L$-theories and elementary classes in language $L$.
\qdr\ect

\bdf\ In an AEC \ $\cK \sub L-str$, define:\\
(i) the "generalized elementary relation" in $\cK$: $\equiv_K$ is "connection relation in $cat(\cK)$", i.e., is the equivalence relation generated by the arrows in $cat(\cK)$:\\
$M \equiv_K N$ iff $\exists n \in \N$, $P_0, \ldots P_n \in \cK$ and $\exists$ zigzag\
 $M= P_0 \overset{f_1}\to\ra P_1 \overset{f_2}\to\la P_2 \overset{f_3}\to\ra  \ldots \overset{f_{n-1}}\to\la P_{n-1} \overset{f_{n}}\to\ra P_n = N$  \\
(ii) a subclass $\cM \sub \cK$ is a "generalized elementary subclass" or a "generalized (closed) theory", when $\cM$ is closed under $\equiv_K$;\\
(iii) there is only a set (of representatives) of closed theories (by LS axiom);\\
(iv) a "maximal consistent g-theory" is an equivalence class module $\equiv_K$ (it is a closed theory) ($Max(\cK) = \cK/\equiv_K$);\\ 
(v) a "complete g-theory" is a non-empty subset of a maximal consistent g-theory. 
\qdr\edf

\bfa\ $\cM \sub \cK$ is g-elementary subclass  \ $\Ra$ \ $(\cM, (\preceq_K)\!\!\rest_M)$ is an AEC, a subAEC of $\cK$,  and 
$LS(\cM) \leqs LS(\cK)$.  
\qdr\efa

\bfa\  Let $\alpha : L \ra L'$, $\alpha^\star\!\rest : \cK' \ra \cK$:\\
(i) $M' \equiv_{K'} N' \ \Ra \ M'^\alpha \equiv_K N'^\alpha$\\
(ii) {\und{Induced Theory}}: Let $\alpha^\star\!\!\rest : \cK' \ra \cK$, $Mod(T') \sub \cK'$, then exists 
$Mod(T')^\alpha$ := the least closed theory of $\cK$  that contains $\{M'^\alpha : M' \in Mod(T')\}$
(iii) a commutative diagram:\\

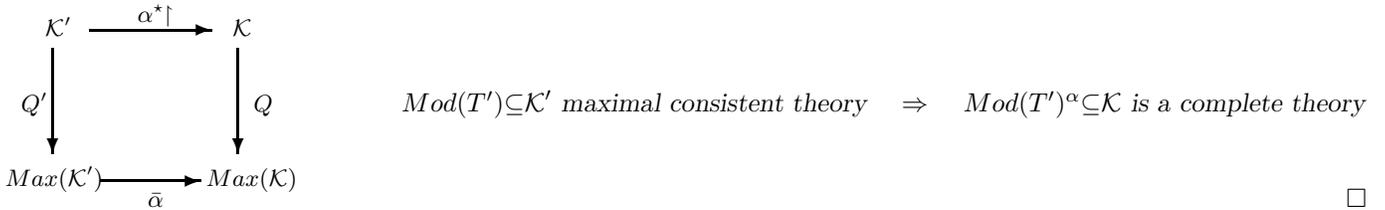
\begin{picture}(100,70)
\setlength{\unitlength}{.4\unitlength} \thicklines 
\put(13,162){\small $\cK'$}
\put(-25,20){\small $Max(\cK')$} 
\put(190,162){\small $\cK$} 
\put(-10,90){\small $Q'$}
\put(100,175){\small $\alpha^\star\!\!\rest$} 
\put(165,20){\small  $Max(\cK)$} 
\put(210,90){\small $Q$} 
\put(110,00){\small $\bar{\alpha}$}
\put(20,150){\vector(0,-1){100}}
\put(195,150){\vector(0,-1){100}} 
\put(55,167){\vector(1,0){115}} 
\put(65,25){\vector(1,0){95}}

\put(350,90){$Mod(T') \sub \cK'$ maximal consistent theory \hem $\Ra$ \hem  $Mod(T')^\alpha \sub \cK$ is a complete theory}
\end{picture}
\qdr\efa

{\bf Limits of AECs}

We will show that the category of all AECs has all finite limits because it has terminal object and pulbacks.  

\bct {\em The terminal AEC:} is obtained taking:\\
\bu $L = \empty$, $\cK = L-str = Set$ \\
\bu $M \preceq_K N$ \ $\Lra$ \ $M \sub N$
\qdr\ect

\bct {\em Pullbacks of AECs:}\\
\bu $L, L_0,L_1$\\
\bu $\alpha_i : L \ra L_i$,  $\alpha_i^\star\rest : (\cK, \preceq_K) \ra (\cK_i, \preceq_{K_i})$ \\
\bu $L' = (L_0 \sqcup L_1)/\sim$ \\
$\pi_i : L_i \ra L'$ (= $L_i \rat L_0 \sqcup L_1 \thra L'$) \\
$\delta = \pi_0 \circ \alpha_0 = \pi_1 \circ \alpha_1 : L \ra L'$\\
$\cK' \sub L'-str$\\
-- $M' \in \cK'$ \ $\Lra$ \
 $M'^{\pi_0} \in \cK_0, M'^{\pi_1} \in \cK_1$, ($(M'^{\pi_0})^{\alpha_0} = (M'^{\pi_1})^{\alpha_1} \in \cK$)\\
-- $M' \preceq_{K'} N'$ \ $\Lra$\
 $M'^{\pi_0} \preceq_{K_0} N'^{\pi_0}, M'^{\pi_1} \preceq_{K_1} N'^{\pi_1}$, ($M'^\delta \preceq_K N'^\delta$)
 
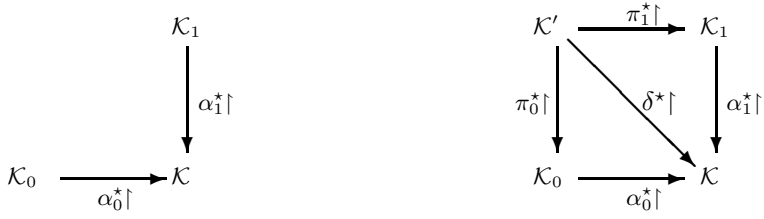
\begin{picture}(100,80)
\setlength{\unitlength}{.4\unitlength} \thicklines 
\put(15,20){\small $\cK_0$} 
\put(170,162){\small $\cK_1$} 
\put(170,20){\small  $\cK$} 
\put(195,90){\small $\alpha_1^\star\!\!\rest$} 
\put(100,0){\small $\alpha_0^\star\!\!\rest$}

\put(185,150){\vector(0,-1){100}} 
\put(65,25){\vector(1,0){100}}

\put(512,162){\small $\cK'$}
\put(512,20){\small $\cK_0$} 
\put(670,162){\small $\cK_1$} 
\put(495,90){\small $\pi_0^\star\!\!\rest$}
\put(600,175){\small $\pi_1^\star\!\!\rest$} 
\put(670,20){\small  $\cK$} 
\put(695,90){\small $\alpha_1^\star\!\!\rest$} 
\put(600,0){\small $\alpha_0^\star\!\!\rest$}
\put(615,90){$\delta^\star\!\!\rest$}
\put(535,150){\vector(0,-1){100}}
\put(685,150){\vector(0,-1){100}} 
\put(555,167){\vector(1,0){100}} 
\put(545,157){\vector(1,-1){120}} 
\put(555,25){\vector(1,0){100}}
\end{picture}

\bu {\em Remark:} $cat(\cK)_{strict}$ is the pullback in $CAT$, the category of all (large) categories.

\bu {{\bf Proposition:}} $\cK'$ is an AEC: \\
\hp \bdm The validity of the axioms {\em order, isomorphism, coherence, reunion} follows from the relation $L'-str \cong (L_0-str)  \underset{L-str}\to\times (L_1-str)$. 
For {\em L\"owenhein-Skolem}: 
$LS(\cK), LS(\cK_0), LS(\cK_1) \leqs \lambda \ \Ra \ LS(\cK') \leqs \lambda $. Take $X \sub M' \in \cK'$ and consider the diagram:
 
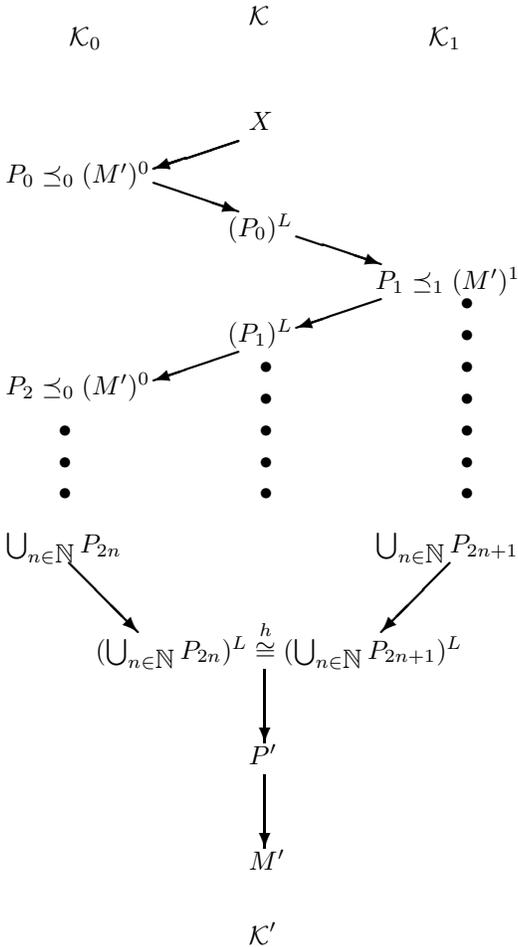
\begin{picture}(300,450)
\setlength{\unitlength}{.4\unitlength} \thicklines 
\put(230,1000){$X$}
\put(220,990){\vector(-3,-1){80}}
\put(220,790){\vector(-3,-1){80}}
\put(355,840){\vector(-3,-1){80}}
\put(140,950){\vector(3,-1){80}}
\put(275,900){\vector(3,-1){80}}
\put(60,590){\vector(1,-1){65}}
\put(420,590){\vector(-1,-1){65}}
\put(245,390){\vector(0,-1){70}}
\put(245,490){\vector(0,-1){70}}
\put(50,710){\bu}
\put(50,680){\bu}
\put(50,650){\bu}
\put(240,710){\bu}
\put(240,680){\bu}
\put(240,650){\bu}
\put(240,740){\bu}
\put(240,770){\bu}
\put(430,740){\bu}
\put(430,710){\bu}
\put(430,680){\bu}
\put(430,650){\bu}
\put(430,770){\bu}
\put(430,800){\bu}
\put(430,830){\bu}
\put(220,990){\vector(-3,-1){80}}
\put(60,1080){$\cK_0$}
\put(400,1080){$\cK_1$}
\put(230,1100){$\cK$}
\put(230,230){$\cK'$}
\put(00,950){$P_0 \preceq_0 (M')^0$}
\put(210,900){$(P_0)^L$}
\put(350,850){$P_1 \preceq_1 (M')^1$}
\put(210,800){$(P_1)^L$}
\put(00,750){$P_2 \preceq_0 (M')^0$}
\put(00,600){$\bigcup_{n \in \N}P_{2n}$}
\put(350,600){$\bigcup_{n \in \N}P_{2n+1}$}
\put(85,500){$(\bigcup_{n \in \N}P_{2n})^{L} \overset{h}\to\cong (\bigcup_{n \in \N}P_{2n+1})^L$}
\put(230,400){$P'$}
\put(230,300){$M'$}

\put(550,770){$P' \in L'-str$ (well defined):}
\put(550,740){$(P')^0 =  \bigcup_{n \in \N}P_{2n} \preceq_0 (M')^0$}
\put(550,710){$(P')^1 \overset{h}\to\cong \bigcup_{n \in \N}P_{2n+1} \preceq_1 (M')^1$}
\put(550,680){$(P')^1 \sub_{L_1} (M')^1$}
\put(550,650){iso condition in $\cK_1$ entails $(P')^1 \preceq_1 (M')^1$}
\put(550,620){thus $P' \preceq_{K'} M'$, $X \sub P'$, $\sharp P' \leqs \lambda$}

\end{picture}

\vspace{-3.0cm}

\bu {\em Special case}:\\
 $L = L_0 \cap L_1$, $L' = L_0 \cup L_1$ 
 (pullback/pushout)

\begin{picture}(100,80)
\setlength{\unitlength}{.4\unitlength} \thicklines 
\put(122,162){\small $(L_0 \cup L_1)-str$}
\put(152,20){\small $L_0-str$} 
\put(570,162){\small $L_1-str$} 
\put(165,90){\small $\pi_0^\star\!\!\rest$}
\put(400,175){\small $\pi_1^\star\!\!\rest$} 
\put(540,20){\small  $(L_0 \cap L_1)-str$} 
\put(630,90){\small $\alpha_1^\star\!\!\rest$} 
\put(400,0){\small $\alpha_0^\star\!\!\rest$}
\put(300,90){}
\put(205,150){\vector(0,-1){100}}
\put(625,150){\vector(0,-1){100}} 
\put(285,167){\vector(1,0){270}} 
\put(255,25){\vector(1,0){270}}
\end{picture}

\qdr\ect

{\bf AECs for g-types}

Similar to pullback (smaller class, but with pullback relation).

\bct
Given $M \in \cK$, consider $\cK_{M}$ the class of all  pairs $(N,f)$ where $N \in \cK$, $ f : M \rat N \in cat(\cK)$.

\bu $\cK_{M} \sub L(M)-str$

\bu binary relation: $(N_0,f_0) \preceq_{K_M} (N_1,f_1)$ iff $N_0 \preceq_K N_1$ and $(N_0,f_0) \sub_{L(M)} (N_1,f_1)$

\bu $(\cK_{M}, \preceq_{K_M})$ is AEC\hem
 ($LS(\cK) \leqs \lambda \Ra LS(\cK_{M}) \leqs \lambda + \sharp M$) 

\bu $cat(\cK_{M})$: commutative triangles in $cat(\cK)$ over $M$

\begin{picture}(100,100)
\thicklines \setlength{\unitlength}{.5\unitlength} \put(5,162){$N_0$}
\put(47,167){\vector(1,0){95}} \put(150,162){$N_1$} \put(20,90){$f_0$}
\put(145,90){$f_1$} \put(80,20){$M$} \put(90,175){$g$}
\put(80,45){\vector(-1,2){55}} \put(105,45){\vector(1,2){55}}
\end{picture} 

\qdr\ect

 \vspace{-1cm}
 
\bct {\em An AEC for g-types}

Given $M \in \cK$ and $I$ a set, consider $\cK_{M,I}$ the class of all  triples $(N,f, \bar{a})$ where $N \in \cK$, $\bar{a} :I \ra N$, $ f : M \rat N \in cat(\cK)$.

\bu $\cK_{M,I} \sub L(M \cup I)-str$

\bu binary relation: $(N_0,f_0, \bar{a}_0) \preceq_{K_{MI}} (N_1,f_1, \bar{a}_1) $ iff $N_0 \preceq_K N_1$ and $(N_0,f_0, \bar{a}_0) \sub_{L({M,I})} (N_1,f_1, \bar{a}_1) $

\bu $(\cK_{M,I}, \preceq_{K_{M,I}})$ is AEC \hem
($LS(\cK) \leqs \lambda \Ra LS(\cK_{M,I}) \leqs \lambda + \sharp M +\sharp I$) 

\bu $cat(\cK_{M,I})$: commutative triangles in $cat(\cK)$ over $M$ and $I$

\bu we may assume $M \cap I = \empty$ (e.g., take $I := C \times \{M\}$)

\begin{picture}(100,100)
\thicklines \setlength{\unitlength}{.5\unitlength} \put(-15,162){$(N_0,\bar{a}_0)$}
\put(77,167){\vector(1,0){95}} \put(180,162){$(N_1,\bar{a}_1)$} \put(50,90){$f_0$}
\put(175,90){$f_1$} \put(110,20){$M$} \put(120,175){$g$}
\put(110,45){\vector(-1,2){50}} \put(135,45){\vector(1,2){50}}
\end{picture} 
\qdr\ect

 \vspace{-1cm}

\bfa The $I$-types over $M$ are precisely the g-maximal consistent theories in $\cK_{M,I}$, thus they are AECs in the language $L(M\cup I)$.\\
\hfl $g-T(M,I)= Max(\cK_{M,I})$\hfl
\qdr\efa

\bfa {\em Induced arrows  on g-types:} \\
$\alpha : L \ra L'$, $\alpha^\star\!\!\rest : \cK' \ra \cK$, $h' : M'_0 \ra M'_1 \in cat(\cK')$, $C$ a set of constants

\bu a commutative diagram of languages and morphisms: \\
-- $I_j = C \times \{|M'_j|\}$, $j =0,1$;\\
-- $M_j := M_j'^\alpha$ $\Ra$ $|(M_j)^\alpha| = |M'_j|$, $j =0,1$;\\
-- $h :=(h')^\alpha : M_0 \ra M_1$ 
 
\begin{picture}(100,80)
\setlength{\unitlength}{.4\unitlength} \thicklines 
\put(132,162){\small $L(M_0 \cup I_0)$}
\put(132,20){\small $L'(M'_0 \cup I_0)$} 
\put(570,162){\small $L(M_1 \cup I_1)$} 
\put(165,90){\small $\td{\alpha}^0$}
\put(400,175){\small $\td{h}$} 
\put(570,20){\small  $L'(M'_1 \cup I_1)$} 
\put(630,90){\small $\td{\alpha}^1$} 
\put(400,-10){\small $\td{h'}$}
\put(300,90){}
\put(205,150){\vector(0,-1){100}}
\put(625,150){\vector(0,-1){100}} 
\put(250,167){\vector(1,0){300}} 
\put(250,25){\vector(1,0){300}}
\end{picture}

\bu a commutative diagram of AECs and morphisms: 
 
\begin{picture}(100,80)
\setlength{\unitlength}{.4\unitlength} \thicklines 
\put(172,162){ $\cK_{M_0, I_0}$}
\put(172,20){ $\cK'_{M'_0, I_0}$} 
\put(570,162){ $\cK_{M_1,I_1}$} 
\put(145,90){\small $(\td{\alpha}^0)^\star$}
\put(400,175){\small $(\td{h})^\star$} 
\put(570,20){  $\cK'_{M'_1, I_1}$} 
\put(630,90){\small $(\td{\alpha}^1)^\star$} 
\put(400,-10){\small $(\td{h'})^\star$}
\put(300,90){}
\put(205,40){\vector(0,1){100}}
\put(615,40){\vector(0,1){100}} 
\put(555,167){\vector(-1,0){300}} 
\put(555,25){\vector(-1,0){300}}
\end{picture}

\bu a commutative diagram of sets and functions:

\begin{picture}(100,80)
\setlength{\unitlength}{.4\unitlength} \thicklines 
\put(132,162){ $Max(\cK_{M_0, I_0})$}
\put(132,20){ $Max(\cK'_{M'_0, I_0})$} 
\put(540,162){ $Max(\cK_{M_1,I_1})$} 
\put(165,90){\small $\bar{\td{\alpha}^0}$}
\put(400,175){\small $\bar{\td{h}}$} 
\put(540,20){  $Max(\cK'_{M'_1, I_1})$} 
\put(630,90){\small $\bar{\td{\alpha}^1}$} 
\put(400,-10){\small $\bar{\td{h'}}$}
\put(300,90){}
\put(205,45){\vector(0,1){100}}
\put(615,45){\vector(0,1){100}} 
\put(530,167){\vector(-1,0){230}} 
\put(530,25){\vector(-1,0){230}}
\end{picture}

 equivalently 

\begin{picture}(100,80)
\setlength{\unitlength}{.4\unitlength} \thicklines 
\put(132,162){ $g-T(M_0, I_0)$}
\put(132,20){ $g-T(M'_0, I_0)$} 
\put(540,162){ $g-T(M_1,I_1)$} 
\put(165,90){\small $\bar{\td{\alpha}^0}$}
\put(400,175){\small $\bar{\td{h}}$} 
\put(540,20){  $g-T(M'_1, I_1)$} 
\put(630,90){\small $\bar{\td{\alpha}^1}$} 
\put(400,-10){\small $\bar{\td{h'}}$}
\put(300,90){}
\put(205,45){\vector(0,1){100}}
\put(615,45){\vector(0,1){100}} 
\put(530,167){\vector(-1,0){230}} 
\put(530,25){\vector(-1,0){230}}
\end{picture}
\qdr\efa

\bfa {\em $\lambda$-Topology on g-types (II)}: 

(i) $\lambda \geqs LS(\cK)$ $\Ra$ $X^\lambda(M,I)$ is the least topology such that:\\
-- $\bar{\td{i_N}} : g-T(M,I) \ra g-T(N,I)$ is continuous, $N \in 
S^\lambda(M) := \{ N \in \cK: N \preceq_K M, \sharp N \leqs \lambda\}$, $X^\lambda(N,I) = (g-T(N,I), discrete)$\\
equivalently:\hfl\\
-- $can^\lambda_M : X^\lambda(M,I) \ra lim_{N \in S^\lambda(M)}  X^\lambda(N,I)$ (with induced pro-topology) is continuous;\\
equivalently:\hfl\\
-- if $M \cong colim_{j \in J} P_j$ (directed colimit), $\sharp P_j \leqs \lambda$\\
$can : X^\lambda(M,I) \ra lim_{j \in J}  X^\lambda(P_j,I) \ \approx $ \
$lim_{N \in S^\lambda(M)}  X^\lambda(N,I)$ (with induced pro-topology) is continuous

(ii) Properties of $\cK$: pro-topology ?, injective?, surjective?\\
--$\lambda$-Tame $\lrsa$ $\all M,I$, $can$ is injective ($T_2$)\\
--"$\lambda$-compact" $\lrsa$ $\all M,I$, $can$ is surjective

(iii) {\em Induced arrows  on g-types:} \\
$\alpha : L \ra L'$, $\alpha^\star\!\!\rest : \cK' \ra \cK$, $h' : M'_0 \ra M'_1 \in cat(\cK')$, $C$ a set of constants,
$\lambda \geqs LS(\cK), LS(\cK')$ 

\bu $M' \cong colim_{j \in J} P'_j$ $\Ra$ $M'^\alpha \cong  colim_{j \in J} {P'_j}^\alpha$\\
$\hat{h}': S^\lambda(M'_0) \ra S^\lambda(M'_1)$

\bu two diagrams of spaces and continuous functions

\begin{picture}(100,80)
\setlength{\unitlength}{.4\unitlength} \thicklines 
\put(102,162){\small $\underset{N_0 \in S^\lambda(M_0)}\to{lim}\!\!\! X^\lambda(N_0, I_0)$}
\put(102,20){\small $\underset{N'_0 \in S^\lambda(M'_0)}\to{lim}\!\!\! X^\lambda(N'_0, I_0)$} 
\put(530,162){\small $\underset{N_1 \in S^\lambda(M_1)}\to{lim}\!\!\! X^\lambda(N_1, I_1)$} 
\put(165,90){\small $\bar{\td{\alpha}^0}$}
\put(400,175){\small $\bar{\td{h}}$} 
\put(530,20){\small $\underset{N'_1 \in S^\lambda(M'_1)}\to{lim}\!\!\! X^\lambda(N'_1, I_1)$} 
\put(640,90){\small $\bar{\td{\alpha}^1}$} 
\put(400,-10){\small $\bar{\td{h'}}$}
\put(300,90){}
\put(205,45){\vector(0,1){100}}
\put(625,45){\vector(0,1){100}} 
\put(530,167){\vector(-1,0){220}} 
\put(530,25){\vector(-1,0){220}}
\end{picture}

 then

\begin{picture}(100,80)
\setlength{\unitlength}{.4\unitlength} \thicklines 
\put(132,162){ $X^\lambda(M_0, I_0)$}
\put(132,20){ $X^\lambda(M'_0, I_0)$} 
\put(550,162){ $X^\lambda(M_1,I_1)$} 
\put(165,90){\small $\bar{\td{\alpha}^0}$}
\put(400,175){\small $\bar{\td{h}}$} 
\put(550,20){  $X^\lambda(M'_1, I_1)$} 
\put(640,90){\small $\bar{\td{\alpha}^1}$} 
\put(400,-10){\small $\bar{\td{h'}}$}
\put(300,90){}
\put(205,45){\vector(0,1){100}}
\put(625,45){\vector(0,1){100}} 
\put(540,167){\vector(-1,0){240}} 
\put(540,25){\vector(-1,0){240}}
\end{picture}

\bu If $\cK'$ has {\bf (AP)} $\Ra$ $\bar{\td{h'}} :  X^\lambda(M'_1,I) \thra X^\lambda(M'_0,I)$ and $\bar{\td{h}} :  X^\lambda(M_1,I) \thra X^\lambda(M_0,I)$ 
are surjective.
\qdr\efa

\subsection{New (global) axioms and applications}

\hp{\bf \large "Local" axioms and applications:}\\
Where "local" means "for each AEC fixed"\\
Ex.:\\
{\bf (JEP)}: Joint Embbeding Property;\\
{\bf (AP)}: Amalgamation Property;\\
{\bf (LRP)}: "Local Robinson's Property": $M_0 \equiv_{\cK} M_1$ \ iff $\exists N \in \cK \exists j_0, j_1 \in cat(\cK)$ 
 $j_0 : M_0 \rat N$, $j_1 : M_1 \rat N$
\\

{\bf Robinson's diagram method:} 

\bct If $\cK$ is an AEC in language $L$, then for each set of "new constants" $L' := L \cup \{ c_i: i \in I\}$ it is available the AEC $\cK'$ in the language $L'$ whose objects are pairs $(M, \bar{a})$, where $\bar{a}$ is a $I$-sequence  of elements of $M$ and $(M, \bar{a}) \preceq_{\cK'} (N, \bar{b})$ iff:\hem
 $M \preceq_{\cK} N$ and $\bar{a} = \bar{b}$\hem 
 ($\Lra$ \ $M \preceq_{\cK} N$ and $(M, \bar{a}) \sub_{L'} (N, \bar{b})$)

\bu $\cK'$ it is  the pullback:

\begin{picture}(100,80)
\setlength{\unitlength}{.4\unitlength} \thicklines 
\put(-15,20){\small $L'-str$} 
\put(170,162){\small $\cK$} 
\put(170,20){\small  $L-str$} 
\put(200,90){\small $incl$} 
\put(100,0){\small $j^\star$}

\put(185,150){\vector(0,-1){100}} 
\put(65,25){\vector(1,0){100}}

\put(512,162){\small $\cK'$}
\put(470,20){\small $L'-str$} 
\put(670,162){\small $\cK$} 
\put(475,90){\small $incl'$}
\put(600,175){\small $j^\star\!\!\rest$} 
\put(670,20){\small  $L-str$} 
\put(690,90){\small $incl$} 
\put(600,0){\small $j^\star$}
\put(535,150){\vector(0,-1){100}}
\put(685,150){\vector(0,-1){100}} 
\put(555,167){\vector(1,0){100}} 
\put(555,25){\vector(1,0){100}}
\end{picture}
\qdr\ect

\bpr {\bf (LRP)} $\exists j : M \rat N \in cat(\cK)$ iff $\exists \bar{b}$ such that $(M, \bar{a}) \equiv_{K(M)} (N, \bar{b})$ where $\cK(M) \sub L(M)-str$

\begin{picture}(100,80)
\thicklines \setlength{\unitlength}{.4\unitlength} \put(-15,162){$(M,\bar{a})$}
\put(190,162){$(N,\bar{b})$} \put(50,90){$t$}
\put(380,130){$(M,\bar{a}) \overset{t\rest}\to{\underset{\cong}\to\ra} (M',\bar{c})  \overset{incl}\to{\underset{\preceq}\to\ra} (P,\bar{c})$}
\put(380,60){$(N,\bar{b}) \overset{h\rest}\to{\underset{\cong}\to\ra} (N',\bar{c})  \overset{incl}\to{\underset{\preceq}\to\ra} (P,\bar{c})$}
\put(185,90){$h$} \put(95,20){$(P, \bar{c})$} 
\put(50,145){\vector(1,-2){50}} \put(195,145){\vector(-1,-2){50}}
\end{picture} 

\vspace{1cm}

\begin{picture}(100,80)
\thicklines \setlength{\unitlength}{.4\unitlength} \put(-15,162){$(M',\bar{c})$}
\put(77,167){\vector(1,0){95}} 
\put(97,177){$\sub_{L(M)}$} 
\put(190,162){$(N',\bar{c})$} \put(50,90){$\preceq$}
\put(380,130){$M,\bar{a}) \ \overset{t\rest}\to{\underset{\cong}\to\ra} \ (M',\bar{c}) \ \overset{incl}\to{\underset{\preceq}\to\ra} \ (N',\bar{c})  \overset{(h\rest)^{-1}}\to{\underset{\cong}\to\ra} \ (N,\bar{b})$}
\put(460,60){$(M,\bar{a}) \ \overset{j}\to{\underset{}\to\rat} \ (N,\bar{b})$}
\put(175,90){$\preceq$} \put(95,20){$(P, \bar{c})$} 
\put(50,145){\vector(1,-2){50}} \put(195,145){\vector(-1,-2){50}}
\end{picture} 
\qdr\epr

{\bf The (global) category $AEC$:}

The Grothendieck's gluying of all AECs

\und{objects}: pairs $(M,\cK)$, where $M \in \cK$

\und{arrows}: pairs $(h, \phi) : (M,\cK) \ra (M',\cK')$ \\
$\phi : L \ra L'$, $\Phi = \phi^\star\rest : \cK' \ra \cK$, $h : M \rat \Phi(M') \in cat(\cK)$ 

\und{identities}: $(id_M,id_L) :  (M,\cK) \ra (M,\cK)$ 

\und{composition}: $(M,\cK) \ \overset{(h,\phi)}\to\ra \ (M',\cK') \ \overset{(h',\phi')}\to\ra \ (M'',\cK'')$ = 
$(M,\cK) \ \overset{(\Phi(h')\circ h,\phi'\circ\phi)}\to\ra  \ (M'',\cK'')$

($M \overset{h}\to\rat \Phi(M') \overset{\Phi(h')}\to\rat \Phi(\Phi'(M''))$)
\qdr

\bct {\em New (global) equivalences in $AEC$}:

-- the connection relation: $(N_0,\cK_0) \equiv (N_1,\cK_1)$ 

-- relative equivalence over $M \in \cK$: \\
$(N_0,\cK_0) \equiv_{(M,\cK)} (N_1,\cK_1)$, with $h_i : M \rat  (N_i)^{\alpha_i}$ 

Ex. $I$-types over $M$ :  $(M,\cK) \ra ((N,\bar{a}), \cK(I))$
\qdr\ect

\vspace{0.5cm}

{ \bf \large New (global) axioms in $AEC$ }

{\em Goal:} Determine interesting subcategories of $AEC$ satisfying interesting (global) axioms.

{\bf (GAP)} "Global Amalgamation Property"\\
entails a simpler description of the relative equivalence over $M \in \cK$: \\
$(N_0,\cK_0) \equiv_{(M,\cK)} (N_1,\cK_1)$, with $h_i : M \rat  (N_i)^{\alpha_i}$ \\
iff $\exists \cK' \exists P' \exists  j_i : N_i \rat P'^{\beta_i}$ such that the diagram below commutes

\begin{picture}(100,80)
\setlength{\unitlength}{.4\unitlength} \thicklines 
\put(152,162){\small $(M,\cK)$}
\put(152,20){\small $(N_0,\cK_0)$} 
\put(570,162){\small $(N_1, \cK_1)$} 
\put(115,90){\small $(h_0,\alpha_0)$}
\put(360,175){\small $(h_1, \alpha_1)$} 
\put(570,20){\small  $(P',\cK')$} 
\put(630,90){\small $(j_1,\beta_1)$} 
\put(360,-10){\small $(j_0,\beta_0)$}
\put(300,90){}
\put(205,150){\vector(0,-1){100}}
\put(625,150){\vector(0,-1){100}} 
\put(260,167){\vector(1,0){300}} 
\put(260,25){\vector(1,0){300}}
\end{picture}
\\

{\bf (TRP)} "Transversal Robinson Property": \\
\hp \bu $\alpha_i : L \ra L_i$,  $i =0,1$\\
$\cK, \cK_0, \cK_1$

\bu  $(M_0)^K \equiv_{K} (M_1)^K$ \ $\Ra$\
 $\exists N_0 \in \cK_0$ such that:\\
-- $M_0 \preceq_{K_0} N_0$ in $\cK_0$\\
-- $\exists j: (M_1)^K \rat (N_0)^K$ in $cat(\cK)$


{\bf (GRP)} "Global Robinson Property":\\
\hp \bu $\alpha_i : L \ra L_i$, $\beta_i :L_i \ra (L_0 \cup L_1)/\sim$, $i =0,1$\\
$\cK, \cK_0, \cK_1,$\
$\cK' = pullback$

\bu  $(M_0)^K \equiv_{K} (M_1)^K$ \ $\Ra$\
 $\exists P' \in \cK'$ such that:\\
-- $M_0 \preceq_{K_0} (P')^{K_0}$ in $\cK_0$\\
-- $\exists h: M_1 \rat (P')^{K_1}$ in $cat(\cK_1)$





\bpr (TRP) $\Ra$ (GRP)\\
\hp \bdm Just take zig-zags as in the proof of LS property in the pullback AEC
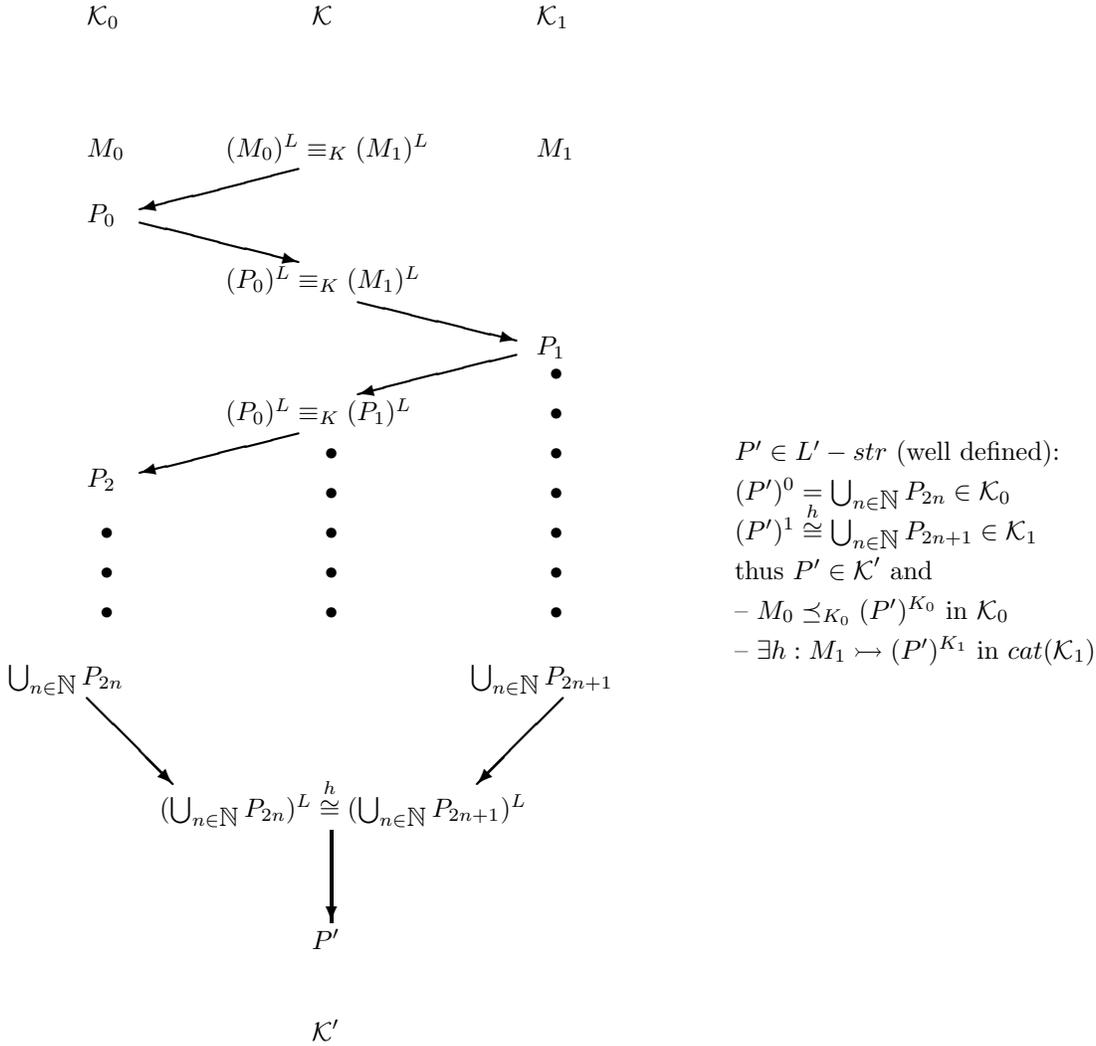
\qdr\epr

\begin{picture}(300,560)
\setlength{\unitlength}{.5\unitlength} \thicklines 
\put(165,1000){$(M_0)^L \equiv_K (M_1)^L$}
\put(60,1000){$M_0$}
\put(400,1000){$M_1$}
\put(220,990){\vector(-4,-1){120}}
\put(220,790){\vector(-4,-1){120}}
\put(385,850){\vector(-4,-1){120}}
\put(100,950){\vector(4,-1){120}}
\put(265,890){\vector(4,-1){120}}
\put(60,590){\vector(1,-1){65}}
\put(420,590){\vector(-1,-1){65}}
\put(245,490){\vector(0,-1){70}}
\put(70,710){\bu}
\put(70,680){\bu}
\put(70,650){\bu}
\put(240,710){\bu}
\put(240,680){\bu}
\put(240,650){\bu}
\put(240,740){\bu}
\put(240,770){\bu}
\put(410,740){\bu}
\put(410,710){\bu}
\put(410,680){\bu}
\put(410,650){\bu}
\put(410,770){\bu}
\put(410,800){\bu}
\put(410,830){\bu}

\put(60,1100){$\cK_0$}
\put(400,1100){$\cK_1$}
\put(230,1100){$\cK$}
\put(230,330){$\cK'$}
\put(60,950){$P_0$}
\put(165,900){$(P_0)^L \equiv_K (M_1)^L$}
\put(400,850){$P_1$}
\put(165,800){$(P_0)^L \equiv_K (P_1)^L$}
\put(60,750){$P_2$}
\put(00,600){$\bigcup_{n \in \N}P_{2n}$}
\put(350,600){$\bigcup_{n \in \N}P_{2n+1}$}
\put(115,500){$(\bigcup_{n \in \N}P_{2n})^{L} \overset{h}\to\cong (\bigcup_{n \in \N}P_{2n+1})^L$}
\put(230,400){$P'$}

\put(550,770){$P' \in L'-str$ (well defined):}
\put(550,740){$(P')^0 =  \bigcup_{n \in \N}P_{2n} \in \cK_0$}
\put(550,710){$(P')^1 \overset{h}\to\cong \bigcup_{n \in \N}P_{2n+1} \in \cK_1$}
\put(550,680){thus $P' \in \cK'$ and}
\put(550,650){-- $M_0 \preceq_{K_0} (P')^{K_0}$ in $\cK_0$}
\put(550,620){-- $\exists h: M_1 \rat (P')^{K_1}$ in $cat(\cK_1)$}

\end{picture}

\vspace{-5.5cm}

{\bf Craig Interpolation Property}:

\bct {\em Craig interpolation in FOL}:

\bu $\psi_0, \psi_1 \in L'-sent$ \
$\psi_0 \vdash \psi_1$ $\Ra$ exists $\theta \in L= L_0 \cap L_1 $: $\psi_0 \vdash \theta$, $\theta \vdash \psi_1$

\bu $\{\psi_0, \neg \psi_1\} \vdash \bot \Ra$ exists $\theta \in L= L_0 \cap L_1 $:  $\psi_0 \vdash \theta$, $\neg \psi_1 \vdash \neg\theta$

\bu $T_0 \cup T_1$ has no $L'$-model   $\Ra$ exists $\theta \in L= L_0 \cap L_1 $: $T_0 \vDash \theta$, $T_1 \vDash \neg\theta$ 

\bu $L = L_0 \cap L_1, T_i\rest_L = \{ \psi \in L-sent : T_i \vdash_{L} \psi\}$\\
$T_0 \cup T_1$ has no $L'$-model $\Ra$\
$ T_0 \rest_L \cup \ T_1\rest_L$ has no $L$-model

\bu $\nexists M' \in L'-str$ such that $(M')^{L_i} \in Mod(T_i)$, $i =0,1$ \ $\Ra$ \
$\nexists  M \in L-str$ with $M \in Mod(T_0\!\!\rest_L) \cap Mod(T_1\!\!\rest_L)$
\qdr\ect

{\bf (CIP)} "Craig Interpolation Property":\\
\hp \bu $\alpha_i : L \ra L_i$, $\beta_i :L_i \ra (L_0 \cup L_1)/\sim$, $i =0,1$\\
$\cK, \cK_0, \cK_1,$\
$\cK' = pullback$\\
$Mod(T_i) \sub \cK_i$, $i=0,1$.

\bu $\nexists M' \in \cK'$ such that $(M')^{K_i} \in Mod(T_i)$, $i =0,1$ \ $\Ra$ \
$\nexists  M \in \cK$ with $M \in Mod(T_0\!\!\rest_K) \cap Mod(T_1\!\!\rest_K)$\
(the induced theories)


\bpr $(TRP) \Ra (CIP)$\\
\hp \bdm (contrapositive)\\
Let $N \in \cK$ with $N \in Mod(T_0\!\!\rest_L) \cap Mod(T_1\!\!\rest_L)$.\\
Then: there are $M_i \in Mod(T_i)$, $i=0,1$, with $(M_0)^K \equiv_K N \equiv_K (M_1)^K$.\\
As $(TRP) \Ra (GRP)$, there is $P' \in \cK'$ such that:\\
\bu $M_0 \preceq_{K_0} (P')^{K_0}$ in $\cK_0$\hfl 
\bu $\exists h: M_1 \rat (P')^{K_1}$ in $cat(\cK_1)$.\hfl\\
Then:\\
\bu $(P')^{K_0} \equiv_{K_0} M_0 \in Mod(T_0)$;\hfl
\bu $(P')^{K_1} \equiv_{K_1} M_1 \in Mod(T_1)$.\hfl
\qdr\epr

\section{Final Remarks and Future Works}

\hp \bu Introduce and study other global axioms (inspirated on model theory of FOL).

\bu Consider functors induced by more general language morphisms $f : L \ra L'$:\\
(i) non injective;\\
(ii) "flexible morphism":\\
  $f_n \in L$ (n-nary functional symbol) $\mapsto$ $t_n \in L'$ (n-ary term)\\
 $R_n \in L$ (n-ary relational symbol) $\mapsto$ $\phi_n \in L'$ (n-ary atomic formula)

\bu Consider more general functors between AECs (ex. functors that are not "over $Set$").

\bu Study metric abstract elementary classes (MAECs) in this global setting.



\end{document}